\documentclass[a4paper,12pt]{article}
\usepackage{amsmath,amssymb,amsfonts}
\makeatletter
\newbox\bk@bxb
\newbox\bk@bxa
\newif\if@bkcont
\newcount\bk@lcnt

\def\breakboxskip{2pt}
\def\breakboxparindent{1.8em}

\def\breakbox{\vskip\breakboxskip\relax
\setbox\bk@bxb\vbox\bgroup
\advance\linewidth -2\fboxrule
\hsize\linewidth\@parboxrestore
\parindent\breakboxparindent\relax}

\def\bk@split{%
\@tempdimb\ht\bk@bxb 
\advance\@tempdimb\dp\bk@bxb
\setbox\bk@bxa\vsplit\bk@bxb to\z@ 
\setbox\bk@bxa\vbox{\unvbox\bk@bxa}
\setbox\@tempboxa\vbox{\copy\bk@bxa\copy\bk@bxb}
\advance\@tempdimb-\ht\@tempboxa
\advance\@tempdimb-\dp\@tempboxa}

\def\bk@addfsepht{%
\setbox\bk@bxa\vbox{\vskip\fboxsep\box\bk@bxa}}

\def\bk@addskipht{%
\setbox\bk@bxa\vbox{\vskip\@tempdimb\box\bk@bxa}}

\def\bk@addfsepdp{%
\@tempdima\dp\bk@bxa
\advance\@tempdima\fboxsep
\dp\bk@bxa\@tempdima}

\def\bk@addskipdp{%
\@tempdima\dp\bk@bxa
\advance\@tempdima\@tempdimb
\dp\bk@bxa\@tempdima}

\def\bk@line{%
\hbox to \linewidth{%
\hskip-2\fboxsep\vrule \@width\fboxrule\hskip.5\fboxsep\vrule \@width\fboxrule\hskip1.5\fboxsep
\box\bk@bxa\hfil
}}%

\def\endbreakbox{\egroup
\ifhmode\par\fi{\noindent\bk@lcnt\@ne
\@bkconttrue\baselineskip\z@\lineskiplimit\z@
\lineskip\z@\vfuzz\maxdimen
\bk@split\bk@addfsepht\bk@addskipdp
\ifvoid\bk@bxb 
\def\bk@fstln{\bk@addfsepdp
\hskip-\parindent\vbox{\llap{\raisebox{-2ex}{\rule{1.5\fboxsep}{\fboxrule}\hskip.5\fboxsep}}\bk@line\llap{\rule{1.5\fboxsep}{\fboxrule}\hskip.5\fboxsep}}}

\else 
\def\bk@fstln{\vbox{\llap{\raisebox{-2ex}{\rule{1.5\fboxsep}{\fboxrule}\hskip.5\fboxsep}}\bk@line}\hfil%
\advance\bk@lcnt\@ne
\loop
\bk@split\bk@addskipdp\leavevmode
\ifvoid\bk@bxb 
\@bkcontfalse\bk@addfsepdp
\vtop{\bk@line\llap{\rule{2\fboxsep}{\fboxrule}}}%

\else 
\bk@line
\fi
\hfil\advance\bk@lcnt\@ne
\if@bkcont\repeat}%
\fi
\leavevmode\bk@fstln\par}\vskip\breakboxskip\relax}

\newcommand{\bizlie}[1]{\mathop{\,\raisebox{-.5ex}{$\widehat{\raisebox{.9ex}{\rule{2.5ex}{.07ex}}}_{#1}$}\,}}

\def\smp{\smallskip\par}
\def\un{{\bf 1}}

\def\pf{\noindent{\bf Proof~:}\ }

\def\findemo{~\leaders\hbox to 1em{\hss\  \hss}\hfill~\raisebox{.5ex}{\framebox[1ex]{}}\smp}

\def\mpn{\medskip\par\noindent}
\def\smpn{\smallskip\par\noindent}
\def\normal{\mathop{\trianglelefteq}}

\def\smp{\smallskip\par}
\def\smpn{\smallskip\par\noindent}

\def\mpoint{\;\;.}
\def\mvirg{\;\;,}

\def\Inf{{\rm Inf}}
\def\Def{{\rm Def}}

\def\Indinf{{\rm Indinf}}

\def\Hom{{\rm Hom}}
\def\End{{\rm End}}

\def\Inf{{\rm Inf}}

\def\Aut{{\rm Aut}}

\def\Id{{\rm Id}}

\def\Irr{{\rm Irr}}

\def\op{^{op}}

\def\N{\mathbb{N}}
\def\F{\mathbb{F}}
\def\Q{\mathbb{Q}}

\newcommand{\dirsum}[1]{\mathop{\oplus}_{#1}\limits}
\newcommand{\romain}[1]{\uppercase\expandafter{\romannumeral #1}}

\newcommand{\sur}[1]{\,\overline{\! #1}}

\newcommand{\prodb}[2]{\mathop{\prod}_{{\scriptstyle #1}\atop{\scriptstyle #2}}}

\def\op{^{op}}

\newenvironment{enonce}[1]{\pagebreak[2]\refstepcounter{subsection}\refstepcounter{prop}\smpn{{\bf \thesection.\arabic{prop}.\ \ #1~:}}\begin{it} }{\end{it}\smp}
\newenvironment{enonce*}[1]{\pagebreak[2]\smpn{#1~:}\begin{it} }{\end{it}\smp}
\newcommand{\result}[1]{\begin{enonce}{#1}}
\def\fresult{\end{enonce}}
\newcommand{\npar}{\smallskip\par\noindent\pagebreak[2]\refstepcounter{subsection}\refstepcounter{prop}{\bf \thesection.\arabic{prop}.\ \ }}

\newenvironment{mth}[1]{\begin{breakbox}\begin{enonce}{#1}}{\end{enonce}\end{breakbox}}
\newenvironment{mth*}[1]{\begin{breakbox}\begin{enonce*}{#1}}{\end{enonce*}\end{breakbox}}
\newenvironment{rem}[1]{\refstepcounter{subsection}\refstepcounter{prop} \mpn{{\bf \thesection.\arabic{prop}.}\ \ \bf#1\ :}}{\smp}

\def\dom{\backslash}
\makeatletter
\renewenvironment{enumerate}{\ifnum \@enumdepth >3 \@toodeep\else
      \advance\@enumdepth \@ne
      \edef\@enumctr{enum\romannumeral\the\@enumdepth}\list
      {\csname label\@enumctr\endcsname}{\setlength{\topsep}{1ex}\setlength{\itemsep}{0pt}\usecounter
        {\@enumctr}\def\makelabel##1{\hss\llap{##1}}}\fi}{\endlist}
\renewenvironment{itemize}{\ifnum \@itemdepth >3 \@toodeep\else \advance\@itemdepth \@ne
\edef\@itemitem{labelitem\romannumeral\the\@itemdepth}%
\list{\csname\@itemitem\endcsname}{\setlength{\topsep}{1ex}\setlength{\itemsep}{0pt}\def\makelabel##1{\hss\llap{##1}}}\fi}
{\endlist}
\def\@sect#1#2#3#4#5#6[#7]#8{\ifnum #2>\c@secnumdepth
    \let\@svsec\@empty\else
    \refstepcounter{#1}\edef\@svsec{\csname the#1\endcsname .\hskip .5em}\fi
    \@tempskipa #5\relax
     \ifdim \@tempskipa>\z@
       \begingroup #6\relax
         \@hangfrom{\hskip #3\relax\@svsec}{\interlinepenalty \@M #8\par}%
       \endgroup
      \csname #1mark\endcsname{#7}\addcontentsline
        {toc}{#1}{\ifnum #2>\c@secnumdepth \else
                     \protect\numberline{\csname the#1\endcsname}\fi
                   #7}\else
       \def\@svsechd{#6\hskip #3\relax  
                  \@svsec #8\csname #1mark\endcsname
                     {#7}\addcontentsline
                          {toc}{#1}{\ifnum #2>\c@secnumdepth \else
                            \protect\numberline{\csname the#1\endcsname}\fi
                      #7}}\fi
    \@xsect{#5}}
\def\section{\@startsection {section}{1}{\z@}{-3.5ex plus-1ex minus
    -.2ex}{2.3ex plus.2ex}{\reset@font\Large\bf}}  

\makeatother
\renewenvironment{equation}{\refstepcounter{subsection}\refstepcounter{prop}$$}{\leqno{\bf (\theprop)}$$}

\def\mar[#1]{\ar@{-}[#1]|-{\object@{<}}}
\def\marb[#1]{\ar@{-}[#1]|{\object+{  }}}

\usepackage[all]{xy}
\def\CB{\mathcal{B}}
\def\GL{{GL}}
\def\Ab{\mathcal{A}b}
\def\Ver{{\rm Ver}}
\newcommand{\bigoplusc}[3]{\mathop{\bigoplus}_{{\scriptstyle #1}\atop {{\scriptstyle #2}\atop {\scriptstyle #3}}}}

\begin{document}
\centerline{\Large\bf $K$-theory, genotypes, and biset functors}\vspace{4ex}
\centerline{\bf Serge Bouc}\vspace{5ex}\par
\begin{center}
\begin{minipage}{9cm}
\begin{footnotesize}
{\bf Abstract: } Let $p$ be an odd prime number. In this paper, we show that the {\em genome} $\Gamma(P)$ of a finite $p$-group $P$, defined as the direct product of the genotypes of all rational irreducible representations of $P$, can be recovered from the first group of $K$-theory $K_1(\Q P)$. It follows that the assignment $P\mapsto \Gamma(P)$ is a $p$-biset functor. We give an explicit formula for the action of bisets on $\Gamma$, in terms of generalized transfers associated to left free bisets. Finally, we show that $\Gamma$ is a rational $p$-biset functor, i.e. that $\Gamma$ factors through the Roquette category of finite $p$-groups.\vspace{2ex}\par
\vspace{2ex}\par
\par
\end{footnotesize}
\end{minipage}
\end{center}
\section{Introduction}
Let $p$ be a prime number. This article originates in a joint work with Nadia Romero (\cite{bouc-romero}), when we started considering the possible applications of {\em genetic bases} to the computation of Whitehead groups of finite $p$-groups. Indeed, after the comprehensive book of B. Oliver (\cite{bob}), it became clear to N. Romero that these questions have close links to rational representations of $p$-groups. So the idea emerged that possibly genetic bases would be a natural tool in this context, and a first use of this is made in~\cite{romero-whitehead}.\par
In particular, when trying to compute various groups related to the Whitehead group of a finite $p$-group $P$ (for odd $p$), a specific product appears, defined in terms of the fields of endomorphisms of the irreducible $\Q P$-modules. After some non trivial reformulation using genetic bases, this product can be viewed as
$$\Gamma(P)=\prod_{S\in\CB}\limits\big(N_P(S)/S\big)\mvirg$$
where $\CB$ is a genetic basis of $P$. As the groups $N_P(S)/S$ are called the {\em types} or {\em genotypes} of the irreducible $\Q P$-modules, we call $\Gamma(P)$ the {\em genome} of $P$. It is the main subject of this paper.\par
The connection of $\Gamma(P)$ with Whitehead groups and $K$-theory is established in Theorem~\ref{genome}: the genome of $P$ can be recovered as the $p$-torsion part of $K_1(\Q P)$. This induces a structure of $p$-biset functor on the correspondence $P\mapsto \Gamma(P)$, which we try to make explicit in Section~\ref{explicit}, by giving formulae to compute the action of a $(Q,P)$-biset on $\Gamma(P)$ (Theorem~\ref{formula}). Finally, we show that $\Gamma$ is a {\em rational} $p$-biset functor, hence it factors through the {\em Roquette category} of finite $p$-groups introduced in \cite{roquette-category}.
\section{Review of $K_1$}
\npar Let $A$ be a ring (with 1). Let $\GL(A)$ denote the colimit of the linear groups $\GL_n(A)$, for $n\in \N_{>0}$, where the inclusion $\GL_n(A)\hookrightarrow \GL_{n+1}(A)$ is 
$$M\in \GL_n(A)\mapsto \left(\begin{array}{cc}M&0\\0&1\end{array}\right)\in \GL_{n+1}(A)\mpoint$$
The group $K_1(A)$ is defined as the abelianization of $\GL(A)$, namely
$$K_1(A)=\GL(A)^{ab}=\GL(A)/[\GL(A),\GL(A)]\mpoint$$
\begin{rem}{Remark}\label{units} In particular there is a canonical group homomorphism from the group $A^\times=\GL_1(A)$ of invertible elements of $A$ to $K_1(A)$, which factors as 
$$\xymatrix{
{A^\times}\ar[r]& *!U(.3){A^\times/[A^\times,A^\times]}\ar[r]^-{\alpha_A}&*!U(.3){K_1(A)}\mpoint
}$$ 
\end{rem}
\npar There is an alternative definition of $K_1(A)$: let $\mathcal{P}(A)$ denote the category of pairs $(P,a)$ of a finitely generated projective (left) $A$-module $P$, and an automorphism $a$ of $P$. A morphism $(P,a)\to (Q,b)$ in $\mathcal{P}(A)$ is a morphism of $A$-modules $f:P\to Q$ such that $b\circ f=f\circ a$.\par
Let $[P,a]$ denote the isomorphism class of $(P,a)$ in $\mathcal{P}(A)$, and let $K_{det}(A)$ denote the Grothendieck group with generators the set of these equivalence classes, and the relations of the following two forms
\begin{itemize} 
\item $[P,a\circ a']=[P,a]+[P,a']$, for any $a,a'\in\Aut_A(P)$,
\item $[Q,b]=[P,a]+[R,c]$ whenever there are morphisms $f:[P,a]\to [Q,b]$ and $g:[Q,b]\to [R,c]$ in $\mathcal{P}(A)$ such that the sequence 
$$0\to P\stackrel{f}{\to}Q\stackrel{g}{\to}R\to 0$$ 
is an exact sequence of $A$-modules (in particular, since $R$ is projective, this sequence splits).
\end{itemize}
If $n\in\N_{>0}$ and $m\in \GL_n(A)$, one can view $m$ as an automorphism of the free module $A^n$. Let $\lambda(m)=[A^n,m]\in K_{det}(A)$.
\begin{mth}{Theorem} The assignment $m\mapsto \lambda(m)$ induces a group isomorphism $K_1(A)\cong K_{det}(A)$.
\end{mth}
\pf See \cite{curtisreinerII} Theorem 40.6.\findemo
\pagebreak[3]
\npar Let now $A$ and $B$ be two rings, and let $L$ be a $(B,A)$-bimodule which is finitely generated and projective as a left $B$-module. If $P$ is a finitely generated projective $A$-module, then $P$ is a direct summand of some free $A$-module $A^n$, and then $L\otimes_AP$ is a direct summand of $L\otimes_AA^n\cong L^n$ as a left $B$-module. Hence $L\otimes_AP$ is a finitely generated projective left $B$-module. Then the functor $P\mapsto L\otimes_AP$ induces a functor $T_L:\mathcal{P}(A)\to\mathcal{P}(B)$ such that
$$T_L\big((P,a)\big)=(L\otimes_AP,L\otimes_Aa)\mpoint$$
One checks easily that the defining relations of $K_{det}(A)$ are preserved by this functor, hence there is a well-defined induced group homomorphism 
$$t_L:K_{det}(A)\to K_{det}(B)$$
sending the class $[P,a]$ to the class $[L\otimes_AP,L\otimes_Aa]$. This group homomorphism is called the (generalized) {\em transfer} associated to the bimodule~$L$.\par
The properties of the tensor product of bimodules now translate to properties of this transfer homomorphism:
\begin{mth}{Proposition} \label{transfer}Let $A,B,C$ be rings. In the following assertions, assume that the bimodules involved are finitely generated and projective as left modules. Then:
\begin{enumerate}
\item if $L\cong L'$ as $(B,A)$-bimodules, then $t_L=t_{L'}$.
\item if $L$ is the $(A,A)$-bimodule $A$, then $t_L=\Id_{K_{det}(A)}$.
\item if $L\cong L_1\oplus L_2$ as $(B,A)$-bimodules, then $t_L=t_{L_1}+t_{L_2}$.
\item if $L$ is a $(B,A)$-bimodule and $M$ is a $(C,B)$-bimodule, then 
$$t_{M}\circ t_L=t_{M\otimes_BL}\mpoint$$
\end{enumerate}
\end{mth}
It follows in particular from (2) and (4) that if $L$ is a $(B,A)$-bimodule inducing a Morita equivalence from $A$ to $B$, then $t_L$ is an isomorphism (more precisely, if $M$ is an $(A,B)$-bimodule such that $M\otimes_BL\cong A$ and $L\otimes_AM\cong B$ as bimodules, then $t_L$ and $t_M$ are inverse to one another).\par
\npar The group $K_1(A)$ has been determined for a number of rings $A$. In particular:
\begin{mth}{Theorem} \label{division} \begin{enumerate}
\item Let $D$ be a division ring. Then $K_1(D)\cong D^\times/[D^\times,D^\times]$.
\item Let $F$ be a field. Then the determinant homomorphism 
$$m\in\GL_n(F)\to \det(m)\in F^\times$$
induces an isomorphism $K_1(F)\cong F^\times$.
\end{enumerate}
\end{mth}
\pf See \cite{curtisreinerII} Theorem 38.32.\findemo
\begin{mth}{Proposition} \label{transferiso} Let $\F$ be a field and $G$ be a finite group of order prime to the characteristic of $\F$. Let $\Irr_\F(G)$ denote a set of representatives of isomorphism classes of irreducible $\F G$-modules, and for $V\in\Irr_\F(G)$, let $D_V=\End_{\F G}(V)$ denote the skew field of endomorphisms of $V$. \par
Then $V$ is an $(\F G,D_V\op)$-bimodule, where the action of $g\in G$ and $f\in D_V$ on $v\in V$ is given by $g\cdot v\cdot f=gf(v)=f(gv)$. Let $V^*$ denote the $\F$-dual of~$V$, considered as a $(D_V\op,\F G)$-bimodule. \par
Then the map
$$\xymatrix{
\tau:K_1(\F G)\ar[r]^-{\prod_V\limits t_{V^*}}&*!U(0.45){\prod_{V\in\Irr_\F(G)}\limits K_1(D_V\op)}
}$$
is a well defined isomorphism of abelian groups, with inverse
$$\xymatrix{
*!U(0.45){\tau':\prod_{V\in\Irr_\F(G)}\limits K_1(D_{V}\op)}\ar[r]^-{\prod_V\limits t_V}&K_1(\F G)\mpoint
}$$
\end{mth}
\pf As $|G|$ is invertible in $\F$, the group algebra $\F G$ is semisimple. Moreover for each $V\in\Irr_\F(V)$, the skew field $D_V\op$ is also a semisimple $\F$-algebra. This shows that $V$ is projective and finitely generated as an $\F G$-module, and that $V^*$ is projective and finitely generated as a $D_V\op$-module (that is $V^*$ is a finite dimensional $D_V$-vector space). Hence the generalized transfer maps $t_V: K_1(\F G)\to K_1(D_V\op)$ and $t_{V^*}:K_1(D_V\op)\to K_1(\F G)$ are well defined.\par
Now for any two finitely generated $\F G$-modules $V$ and $W$, the map
$$\alpha\otimes w\mapsto \big(v\in V\mapsto \alpha(v)w\in W\big)$$
extends to an isomorphism (see e.g. \cite{curtisreiner} (2.32))
$$V^*\otimes_{\F G}W\to \Hom_{\F G}(V,W)$$
of $\big((\End_{\F G}V)\op,(\End_{\F G}W\op)\big)$-bimodules, where the bimodule structure on the right hand side is given by 
$$\forall h\in (\End_{\F G}V)\op,\;\forall \psi\in \Hom_{\F G}(V,W),\;\forall k\in(\End_{\F G}W)\op,\;\;h\cdot\psi \cdot k=k\circ\psi\circ h.$$
In case $V,W\in\Irr_\F(G)$ and $V\neq W$, this yields $V^*\otimes_{\F G}W=0$. And if $V=W$, we have an isomorphism $V^*\otimes_{\F G}V\cong D_{V}\op$ of $(D_V\op,D_V\op)$-bimodules. Then by Assertions (2) and (4) of Proposition~\ref{transfer} 
$$t_{V^*}\circ t_W=\left\{\begin{array}{ll}0&\hbox{if}\;V\neq W\\\Id_{K_1(D_V\op)}&\hbox{if}\;V=W.\end{array}\right.$$
In other words $\tau\circ \tau'$ is the identity map of $\prod_{V\in\Irr_\F(G)}\limits K_1(D_{V}\op)$.
Conversely
$$\tau'\circ \tau=\sum_{V\in \Irr_{\F}(G)}t_V\circ t_{V^*}=t_L\mvirg$$
where $L$ is the $(\F G,\F G)$-bimodule $\dirsum{V\in\Irr_{\F}(G)}(V\otimes_{D_V\op}V^*)$. For each $V\in\Irr_{\F}(G)$, the bimodule $V\otimes_{D_V\op}V^*\cong \End_{D_V\op}(V)$ is isomorphic to the Wedderburn component of $\F G$ corresponding to the simple module $V$, and the semisimple algebra $\F G$ is equal to the direct sum of its Wedderburn components. Thus $L\cong \F G$, and $t_L$ is equal to the identity map of $K_1(\F G)$.\findemo
\begin{mth}{Corollary} Under the assumptions of Proposition~\ref{transferiso}, there is a group isomorphism
$$K_1(\F G)\cong\prod_{V\in\Irr_{\F}(G)}D_V^\times/[D_V^\times,D_V^\times]\mpoint$$
\end{mth}
\pf This follows from Proposition~\ref{transferiso} and Theorem~\ref{division}, since $x\mapsto x^{-1}$ is a group isomorphism $D^\times \to (D\op)^\times$, for any skew field $D$.\findemo
\npar Recall (\cite{bisetfunctors} Chapter 3) that the {\em biset category} $\mathcal{C}$ of finite groups has all finite groups as objects, the set of morphisms in $\mathcal{C}$ from a group $G$ to a group $H$ being the Grothendieck group of (finite) $(H,G)$-bisets, i.e. the Burnside group $B(H,G)$. The composition of morphisms in $\mathcal{C}$ is the linear extension of the product $(V,U)\mapsto V\times_HU$, for a $(K,H)$-biset $V$ and an $(H,G)$-biset~$U$. A {\em biset functor} is an additive functor from $\mathcal{C}$ to the category $\Ab$ of abelian groups.\par
For a prime number $p$, a {\em $p$-biset functor} is an additive functor from the full subcategory $\mathcal{C}_p$ of $\mathcal{C}$ consisting of $p$-groups to $\Ab$.\par
Let $_1\mathcal{C}$ denote the (non full) subcategory of $\mathcal{C}$ with the same objects, but where the set of morphisms from a group $G$ to a group $H$ is the Grothendieck group $_1B(H,G)$ of {\em left free} $(H,G)$-bisets. A {\em deflation biset functor} is an additive functor from $_1\mathcal{C}$ to $\Ab$.
\begin{mth}{Proposition} \label{bisetfunctors} \begin{enumerate}
\item Let $R$ be a commutative ring. The assignment $G\mapsto K_1(RG)$ is a deflation functor.
\item The assignment $G\mapsto K_1(\Q G)$ is a biset functor.
\end{enumerate}
\end{mth}
\pf For Assertion~1, if $G$ and $H$ are finite groups, and if $U$ is a finite left free $(H,G)$-biset, then the corresponding permutation $(RH,RG)$-bimodule $RU$ is free and finitely generated as a left $RH$-module. Hence the transfer $t_{RU}:K_1(RG)\to K_1(RH)$ is well defined. If $U'$ is an $(H,G)$-biset isomorphic to $U$, then $RU'\cong RU$ as bimodules, hence $t_{RU'}=t_{RU}$. And if $U$ is the disjoint unions of two $(H,G)$-bisets $U_1$ and $U_2$, then $RU\cong RU_1\oplus RU_2$, thus $t_{RU}=t_{RU_1}+t_{RU_2}$. This shows that one can extend linearly this transfer construction $U\mapsto t_{RU}$ to a group homomorphism 
$$u\in {_1B}(H,G)\mapsto K_1(u)\in\Hom_{\Ab}\big(K_1(RG),K_1(RH)\big)\mpoint$$
Moreover, if $K$ is a third group, and $V$ is a finite left free $(K,H)$-biset, then $t_{RV}\circ t_{RU}=t_{RV\otimes_{RH} RU}=t_{R(V\times_HU)}$ since the bimodules $RV\otimes_{RH} RU$ and $R(V\times_HU)$ are isomorphic. Finally, if $U$ is the identity biset at $G$, namely the set $G$ acted on by left and right multiplication, then $RU\cong RG$ as $(RG,RG)$-bimodule, thus $t_{RU}=\Id_{K_1(RG)}$. This completes the proof of Assertion (1).\par
The proof of Assertion (2) is the same, except that the transfer $t_{\Q U}:K_1(\Q G)\to K_1(\Q H)$ is well defined for an arbitrary finite $(H,G)$-biset $U$: indeed $\Q U$ is always finitely generated and projective as a $\Q H$-module.\findemo 
\vspace{-3ex}
\section{Review of genetic subgroups}
\vspace{-1ex}
\npar Let $p$ be a prime number. A finite $p$-group is called a {\em Roquette} $p$-group if it has normal rank 1, i.e. if all its normal abelian subgroups are cyclic. The Roquette $p$-groups (see~\cite{roquette}) are the cyclic groups $C_{p^n}$, for $n\in \N$, if $p$ is odd. The Roquette 2-groups are the cyclic groups $C_{2^n}$, for $n\in \N$, the generalized quaternion groups $Q_{2^n}$, for $n\geq 3$, the dihedral groups $D_{2^n}$, for $n\geq 4$, and the semidihedral groups $SD_{2^n}$, for $n\geq 4$.\par
If $P$ is a Roquette $p$-group, then $P$ admits a {\em unique faithful irreducible rational representation} $\Phi_P$ (\cite{bisetfunctors} Proposition 9.3.5).
\npar If $S$ is a subgroup of a finite $p$-group $P$, denote by $Z_P(S)$ the subgroup of $N_P(S)$ defined by $Z_P(S)/S=Z\big(N_P(S)/S\big)$. The subgroup $S$ is called {\em genetic} if it fulfills the following two conditions:
\begin{enumerate}
\item if $x\in P$, then $S^x\cap Z_P(S)\leq S$ if and only if $S^x=S$.
\item the group $N_P(S)/S$ is a Roquette $p$-group.
\end{enumerate}
When $S$ is a genetic subgroup of $P$, let $V(S)=\Indinf_{N_P(S)/S}^P\Phi_{N_P(S)/S}$ denote the $\Q P$-module obtained by inflation of $\Phi_{N_P(S)/S}$ to $N_P(S)$ followed by induction to $P$.\par
Two genetic subgroups $S$ and $T$ of $P$ are said to be {\em linked modulo $P$}  (notation $S\bizlie{P}T$) if there exists an element $x\in P$ such that $S^x\cap Z_P(T)\leq T$ and ${^xT}\cap Z_P(S)\leq S$ (where as usual $S^x=x^{-1}Sx$ and ${^xT}=xTx^{-1})$. 
\begin{mth}{Theorem}\label{roquetteplus}  Let $p$ be a prime number and $P$ be a finite $p$-group.
\begin{enumerate}
\item If $V$ is a simple $\Q P$-module, then there exists a genetic subgroup $S$ of $P$ such that $V\cong V(S)$.
\item If $S$ is a genetic subgroup of $P$, then there is an isomorphism of $\Q$-algebras
$$\End_{\Q P}V(S)\cong\End_{\Q N_P(S)/S}\Phi_{N_P(S)/S}$$
induced by the induction-inflation functor from $\Q N_P(S)/S$-modules to $\Q P$-modules.
\item If $S$ and $T$ are genetic subgroups of $P$, then $V(S)\cong V(T)$ if and only if $S\bizlie{P}T$. In this case, the groups $N_P(S)/S$ and $N_P(T)/T$ are isomorphic.
\end{enumerate}
\end{mth}
\pf See Theorem 9.4.1, Lemma 9.4.3, Definition 9.4.4, Corollary 9.4.5, Theorem 9.5.6 and Theorem 9.6.1 of \cite{bisetfunctors}.\findemo
It follows in particular that the relation $\bizlie{P}$ is an equivalence relation on the set of genetic subgroups of $P$. A {\em genetic basis} of $P$ is by definition a set of representatives of genetic subgroups of $P$ for this equivalence.\par
It also follows that if $V$ is a simple $\Q P$-module, and if $S$ is a genetic subgroup of $P$ such that $V\cong V(S)$, then the group $N_P(S)/S$ does not depend on the choice of such a genetic subgroup $S$. This factor group is called the {\em type} of $V$ (\cite{bisetfunctors} Definition 9.6.8). Laurence Barker (\cite{barker-genotype}) has introduced the word {\em genotype} instead of type, and we will follow this terminology.
\begin{mth}{Definition} Let $p$ be a prime number and $P$ be a finite $p$-group. The {\em genome} $\Gamma(P)$ of $P$ is the product group 
$$\Gamma(P)=\prod_{S\in\mathcal{B}}\limits\big(N_P(S)/S\big)\mvirg$$
where $\mathcal{B}$ is a genetic basis of $P$. It is well defined up to isomorphism.
\end{mth}
More precisely, suppose that $\mathcal{B}$ and $\mathcal{B}'$ are genetic bases of a $p$-group $P$. Then for $S\in \mathcal{B}$, there exists a unique $S'\in \mathcal{B}'$ such that there exists some $x\in P$ with
\begin{equation}\label{link}
S^x\cap Z_P(S')\leq S'\;\hbox{and}\;{^xS'}\cap Z_P(S)\leq S\mvirg
\end{equation}
and the correspondence $S\mapsto S'$ is a bijection from $\mathcal{B}$ to $\mathcal{B}'$. Moreover, for each $S\in\mathcal{B}$ corresponding to $S'\in\mathcal{B}'$, the set $\mathcal{D}$ of elements $x$ satisfying~(\ref{link}) is a single $\big(N_P(S),N_P(S')\big)$-double coset in $P$ (\cite{bisetfunctors}, Proposition 9.6.9). \par
Let $x\in\mathcal{D}$. Then for each $n\in N_P(S)/S$, there is a unique element $n'\in N_P(S')/S'$ such that $nSx=xS'n'$, and the map $n\mapsto n'$ is a group isomorphism $N_P(S)/S\to N_P(S')/S'$, which only depends on $x$ up to interior automorphism of $N_P(S)/S$. In particular, when $p$ is odd, the group $N_P(S)/S$ is cyclic, so this group isomorphism does not depend on $x$.\par
Thus for odd $p$, this yields a canonical group isomorphism
\begin{equation}\label{canonical}
\xymatrix{
*!U(0.3){\prod_{S\in\mathcal{B}}\limits\big(N_P(S)/S\big)}\ar[r]^-{\gamma_{\mathcal{B}',\mathcal{B}}}&*!U(0.3){\prod_{S'\in\mathcal{B}'}\limits\big(N_P(S')/S'\big)\mpoint}
}
\end{equation}
\begin{rem}{Remark} Let $p$ be a prime number, and $P$ be a finite $p$-group. Since the Roquette $p$-groups are all indecomposable (that is, they cannot be written as a direct product of two non-trivial of their subgroups), the genotypes of the simple $\Q P$-modules are determined by the group $\Gamma(P)$: by the Krull-Remak-Schmidt theorem, the group $\Gamma(P)$ can be written as a direct product of indecomposable groups $\Gamma_1,\ldots,\Gamma_r$, and such a decomposition is unique (up to permutation and isomorphism of the factors). Then $\Gamma_1,\ldots,\Gamma_r$ are the genotypes of the simple $\Q P$-modules.\par
In terms of the Roquette category $\mathcal{R}_p$ (see Section~\ref{roquette}, or~\cite{roquette-category}), this means that two finite $p$-groups $P$ and $Q$ become isomorphic in $\mathcal{R}_p$ if and only if their genomes $\Gamma(P)$ and $\Gamma(Q)$ are isomorphic (as groups) (see \cite{roquette-category} Proposition~5.14).
\end{rem}
\section{$K$-theory and genome}
\vspace{-2ex}
\begin{mth}{Lemma}\label{Phi cyclic} Let $p$ be a prime, and $C$ be a cyclic $p$-group. Recall that $\Phi_C$ is the unique faithful irreducible rational representation of $C$, up to isomorphism.
\begin{enumerate}
\item If $C=\un$, then $\Phi_C=\Q$.
\item If $C\neq \un$, let $Z$ be the unique subgroup of order $p$ of $C$. Then
there is an exact sequence
\begin{equation}\label{defPhi}
0\to \Phi_C\to \Q C\to \Q (C/Z)\to 0\mvirg
\end{equation}
of $(\Q C,\Q C)$-bimodules, where $\Q C\to \Q (C/Z)$ is the canonical surjection.
\item If $C$ has order $p^n$, then the algebra $\End_{\Q C}(\Phi_C)$ is isomorphic to the cyclotomic field $\Q(\zeta_{p^n})$, and if $p>2$, the map sending $c\in C$ to the endomorphism $\varphi\mapsto\varphi c$ of $\Phi_C$ is a group isomorphism from $C$ to the $p$-torsion part $_p\Q(\zeta_{p^n})^\times$ of the multiplicative group $\Q(\zeta_{p^n})^\times$.
\end{enumerate}
\end{mth}
\pf Assertion~1 is trivial. Assertion~2 follows e.g. from~\cite{bisetfunctors}, Proposition~9.3.5. A different proof consists in observing that if $C$ has order $p^n$, then the algebra $\Q C$ is isomorphic to $\Q[X]/(X^{p^n}-1)$, and the projection map $\Q C\to \Q(C/Z)$ becomes the canonical map
$$\Q[X]/(X^{p^n}-1)\to \Q[X]/(X^{p^{n-1}}-1)\mpoint$$
The kernel of this map is now clearly isomorphic to $\Q[X]/(\gamma_{p^n})$, where $\gamma_{p^n}$ is the $p^n$-th cyclotomic polynomial, that is, the $p^n$-th cyclotomic field, which is clearly a simple faithful module for the cyclic group generated by $X$ in the algebra $\Q[X]/(X^{p^n}-1)$. Observe moreover that the exact sequence~\ref{defPhi} is indeed a sequence of $(\Q C,\Q C)$-bimodules. \par
The first part of Assertion 3 follows easily. For the last part, let $\zeta_{p^n}$ be a primitive $p^n$-th root of unity. Observe that a $p$-torsion element in $\Q(\zeta_{p^n})^\times$ is a $p^n$-th root of unity. Hence the $p$-torsion part of $\Q(\zeta_{p^n})^\times$ is cyclic of order~$p^n$, generated by $\zeta_{p^n}$.\findemo
\begin{mth}{Theorem}\label{genome} Let $p$ be an odd prime, and $P$ be a finite $p$-group, and $\mathcal{B}$ be a genetic basis of $P$. If $S$ is a genetic subgroup of $P$, and $a\in N_P(S)/S$, view $a$ as an automorphism of $\Phi_{N_P(S)/S}$, and let $\tilde{a}$ denote the corresponding automorphism of $V(S)=\Indinf_{N_P(S)}^P\Phi_{N_P(S)/S}$. 
\begin{enumerate}
\item The group homomorphism
$$\xymatrix{
*!U(.4){\Gamma(P)=\prod_{S\in\mathcal{B}}\limits\big(N_P(S)/S\big)}\ar[r]^-{\nu_{\mathcal{B}}}& K_1(\Q P)}$$
sending $a\in N_P(S)/S$, for $s\in\mathcal{B}$, to the class $[V(S),\tilde{a}]$ in $K_1(\Q P)$ is an isomorphism of the genome $\Gamma(P)$ onto the $p$-torsion part $_pK_1(\Q P)$ of $K_1(\Q P)$. 
\item If $\mathcal{B}'$ is another genetic basis of $P$, and $\gamma_{\mathcal{B'},\mathcal{B}}$ is the canonical isomorphism defined in~\ref{canonical}, then
$$\nu_{\mathcal{B}'}\circ \gamma_{\mathcal{B'},\mathcal{B}}=\nu_{\mathcal{B}}\mpoint$$
\end{enumerate}
\end{mth}
\pf Since $p$ is odd, the Roquette $p$-groups are the cyclic $p$-groups. Assertion~1 now follows from Proposition~\ref{transferiso}, Theorem~\ref{roquetteplus}, and Lemma~\ref{Phi cyclic}.\par
For Assertion~2, let $S\in\mathcal{B}$ and let $S'$ be the unique element of $\mathcal{B}'$ such that $S'\bizlie{P}S$. Let $\varphi:N_P(S)/S\to N_P(S')/S'$ be the restriction of $\gamma_{\mathcal{B}',\mathcal{B}}$ to $N_P(S)/S$. If $a\in N_P(S)/S$, let $a'=\varphi(a)$. Then $\varphi$ induces an isomorphism of $\Q P$-modules $\tilde{\varphi}:V(S)\to V(S')$ such that the diagram
$$\xymatrix{
V(S)\ar[d]_-{\tilde{\varphi}}\ar[r]^-{\tilde{a}}&V(S)\ar[d]^-{\tilde{\varphi}}\\
V(S')\ar[r]^-{\tilde{a}'}&V(S')\\
}$$
is commutative. Hence $\big(V(S),\tilde{a}\big)\cong\big(V(S'),\tilde{a}'\big)$ in $\mathcal{P}(\Q P)$, thus $\big[V(S),\tilde{a}\big]=\big[V(S'),\tilde{a}'\big]$ in $K_1(\Q P)$, as was to be shown.  \findemo
\begin{rem}{Remark} The elements of odd order of $\Q(\zeta_{p^n})^\times$ are the $p^n$-th roots of unity. So $\Gamma(P)$ is also the odd-torsion part of $K_1(\Q P)$.
\end{rem}
\begin{mth}{Corollary}\label{genome functor} Let $p$ be an odd prime. Then the correspondence sending a finite $p$-group $P$ to its genome $\Gamma(P)$ is a $p$-biset functor.
\end{mth}
\pf Indeed by Proposition~\ref{bisetfunctors}, the assignment $P\mapsto K_1(\Q P)$ is a $p$-biset functor. So its $p$-torsion part is also a $p$-biset functor.\findemo
\vspace{-3ex}
\section{Explicit transfer maps}\label{explicit}
We begin with a slight generalization of the transfer homomorphism, associated to a left-free biset:
\begin{mth}{Lemma and Definition}\label{verlagerung} Let $G$ and $H$ be finite groups, and let $\Omega$ be a left free $(H,G)$-biset. Let $[H\dom \Omega]$ be a set of representatives of $H$-orbits on~$\Omega$. For $g\in G$, and $x\in\Omega$, let $h_{g,x}\in H$ and $\sigma_g(x)\in[H\dom\Omega]$ be the elements defined by $xg=h_{g,x}\sigma_g(x)$.
\begin{enumerate}
\item The map $g\in G\mapsto \prod_{x\in[H\dom\Omega]}\limits h_{g,x}$ (in any order) induces a well defined group homomorphism
$$\Ver_\Omega:G/[G,G]\to H/[H,H]$$
called the {\em (generalized) transfer} associated to $\Omega$.
\item If $\Omega'\cong\Omega$ as $(H,G)$-bisets, then $\Ver_{\Omega'}=\Ver_{\Omega}$.
\item If $\Omega=\Omega_1\sqcup\Omega_2$ as $(H,G)$-bisets, then $\Ver_{\Omega}=\Ver_{\Omega_1}+\Ver_{\Omega_2}$.
\item If $K$ is another finite group, and $\Omega'$ is a finite left free $(K,H)$-biset, then $\Omega'\times_H\Omega$ is a finite left free $K$-set, and 
$$\Ver_{\Omega'}\circ\Ver_{\Omega}=\Ver_{\Omega'\times_H\Omega}\mpoint$$
\end{enumerate}
\end{mth}
The notation and terminology comes from the classical transfer from $G/[G,G]$ to $H/[H,H]$, when $H$ is a subgroup of $G$: the corresponding biset $\Omega$ is the set $G$ itself, in this case.\vspace{1ex}\par\noindent
\pf Changing the set of representatives $[H\dom \Omega]$ amounts to replacing each $x\in[H\dom\Omega]$ by $\eta_xx$, for some $\eta_x\in H$. This changes the element $h_{g,x}$ in $h'_{g,x}=\eta_xh_{g,x}\eta_{\sigma_g(x)}^{-1}$, so the product over $x\in[H\dom \Omega]$ of the elements $h'_{g,x}$ is equal to the product of the elements $h_{g,x}$ in the abelianization $H/[H,H]$. Hence $\Ver_\Omega$ does not depend on the choice of a set of representatives.\par
It follows moreover from the definition that for $g,g'\in G$ and $x\in [H\dom \Omega]$, we have $h_{gg',x}=h_{g,x}h_{g',\sigma_g(x)}$. Hence
$$\prod_{x\in [H\dom\Omega]}\limits h_{gg',x}=\prod_{x\in [H\dom\Omega]}\limits h_{g,x}\prod_{x\in [H\dom\Omega]}\limits h_{g',\sigma_g(x)}=\prod_{x\in [H\dom\Omega]}\limits h_{g,x}\prod_{x\in [H\dom\Omega]}\limits h_{g',x}$$
in $H/[H,H]$, so $\Ver_\Omega$ is a group homomorphism. This proves Assertion~1.\par
For Assertion 2, let $f:\Omega\to\Omega'$ be an isomorphism of $(H,G)$-bisets. Then the set $f\big([H\dom \Omega]\big)$ is a set of representatives of the $H$-orbits on $\Omega'$. Moreover for $x\in[H\dom \Omega]$ and $g\in G$, 
$$f(x)g=f(xg)=f\big(h_{g,x}\sigma_g(x)\big)=h_{g,x}f\big(\sigma_g(x)\big)\mvirg$$
so $\Ver_{\Omega'}(g)=\prod_{x\in\Omega}\limits h_{g,x}=\Ver_\Omega(g)$, which proves Assertion 2.\par
Assertion 3 is clear, since $[H\dom\Omega]=[H\dom\Omega_1]\sqcup[H\dom\Omega_2]$.\par
For Assertion 4, it is straightforward to check that $\Omega'\times_H\Omega$ is left free. Moreover, the set of pairs $(x',x)\in \Omega'\times_H\Omega$, for $x'\in[K\dom\Omega']$ and $x\in[H\dom\Omega]$, is a set of representatives of $K$ orbits on $\Omega'\times_H\Omega$. Then for $x'\in[K\dom\Omega']$ and $x\in[H\dom\Omega]$, and $g\in G$
\begin{eqnarray*}(x',x)g&=&(x',xg)=\big(x',h_{g,x}\sigma_{g}(x)\big)\\
&=&\big(x'h_{g,x},\sigma_g(x)\big)=\big(k_{h_{g,x},x'}\tau_{h_{g,x}}(x'),\sigma_g(x)\big)\\
&=&k_{h_{g,x},x'}\big(\tau_{h_{g,x}}(x'),\sigma_g(x)\big)\mvirg
\end{eqnarray*}
where $k_{h,x'}\in K$ and $\tau_{h}(x')\in[K\dom\Omega']$ are defined by $x'h=k_{h,x'}\tau_{h}(x')$, for $h\in H$ and $x'\in[K\dom \Omega']$.\par
It follows that
$$\Ver_{\Omega'\times_H\Omega}(g)=\prodb{x\in[H\dom \Omega]}{x'\in[K\dom \Omega']}k_{h_{g,x},x'}=\Ver_{\Omega'}\big(\prod_{x\in[H\dom \Omega]}h_{g,x}\big)=\Ver_{\Omega'}\circ\Ver_{\Omega}(g)\mvirg$$
which completes the proof.  \findemo\vspace{2ex}
\npar\label{functor} Corollary~\ref{genome functor} shows that there exists a $p$-biset functor structure on the assignment $P\mapsto \Gamma(P)$ for $p$-groups, when $p$ is odd. This raises the following question: suppose that $P$ and $Q$ are finite $p$-groups, that $\mathcal{B}_P$ is a genetic basis of $P$, and $\mathcal{B}_Q$ is a genetic basis of $Q$. When $U$ is a finite $(Q,P)$-biset, how can we compute the map
$$\Gamma(U):\Gamma(P)=\prod_{S\in \mathcal{B}_P}\big(N_P(S)/S\big)\to \Gamma(Q)=\prod_{T\in \mathcal{B}_Q}\big(N_Q(T)/T\big)$$
giving the action of the biset $U$? \par
This amounts to finding the map 
$$\Gamma(U)_{T,S}:\sur{N}_P(S)=N_P(S)/S\to \sur{N}_Q(T)=N_Q(T)/T$$
for each pair $(T,S)$ of a genetic subgroup $T$ of $Q$ and a genetic subgroup $S$ of~$P$, defined as follows: if $a\in \sur{N}_P(S)$, then $a$ can be viewed as an automorphism of the $\Q \sur{N}_P(S)$-module $\Phi_{\sur{N}_P(S)}$, viewed as an ideal of $\Q \sur{N}_P(S)$ as in~\ref{defPhi}. Then $\tilde{a}=\Indinf_{\sur{N}_P(S)}^Pa$ is an automorphism of $V(S)=\Indinf_{\sur{N}_P(S)}^P\Phi_{\sur{N}_P(S)}$, hence an element $\hat{a}=[V(S),\tilde{a}]$ of $K_1(\Q P)$. This element is mapped by $t_{\Q U}$ to the element 
$$t_{\Q U}(\hat{a})=[\Q U\otimes_{\Q P}V(S),\Q U\otimes_{\Q P}\tilde{a}]$$ 
of $K_1(\Q Q)$. This in turn is mapped to the element $t_{V(T)^*}\circ t_{\Q U}(\hat{a})$ of the direct summand $K_1(F_T)$ of $K_1(\Q Q)$ corresponding to the simple $\Q Q$-module $V(T)$ as in Proposition~\ref{transferiso}, where $F_T$ is the field $D_{V(T)}=\End_{\Q Q}V(T)$.\par
Thus to find $\Gamma(U)_{T,S}(a)$, we have to compute the element
$$[V(T)^*\otimes_{\Q Q}\Q U\otimes_{\Q P}V(S), V(T)^*\otimes_{\Q Q}\Q U\otimes_{\Q P}\tilde{a}]$$
of $K_1(F_T)\cong F_T^\times$, and identify it as an element of $\sur{N}_Q(T)$.\vspace{1ex}\par
We set $L(U)_{T,S}=V(T)^*\otimes_{\Q Q}\Q U\otimes_{\Q P}V(S)$ for simplicity. First we observe that the induction-inflation functor $\Indinf_{\sur{N}_P(S)}^P$ is isomorphic to the functor $\Q (P/S)\otimes_{\sur{N}_P(S)}({-})$, where $\Q (P/S)$ is endowed with its natural structure of $\Big(\Q P,\Q \big(\sur{N}_P(S)\big)\Big)$-bimodule. Hence
\begin{eqnarray*}
\Q U\otimes_{\Q P}V(S)&=&\Q U\otimes_{\Q P}\Indinf_{\sur{N}_P(S)}^P\Phi_{\sur{N}_P(S)}\\
&\cong& \Q U\otimes_{\Q P}\Q(P/S)\otimes_{\Q \sur{N}_P(S)}\Phi_{\sur{N}_P(S)}\\
&\cong& \Q (U/S)\otimes_{\Q \sur{N}_P(S)}\Phi_{\sur{N}_P(S)}\mvirg
\end{eqnarray*}
where $\Q (U/S)$ is given its natural structure of $\Big(\Q P,\Q \big(\sur{N}_P(S)\big)\Big)$-bimodule.\par
Tensoring on the left with $V(T)^*$, an using a similar argument, we get that
$$L(U)_{T,S}\cong \Phi_{\sur{N}_Q(T)}^*\otimes_{\Q \sur{N}_Q(T)}\Q (T\dom U/S)\otimes_{\Q \sur{N}_P(S)}\Phi_{\sur{N}_P(S)}\mvirg$$
where $\Q (T\dom U/S)$ is the permutation $\big(\Q\sur{N}_Q(T),\Q\sur{N}_P(S)\big)$-bimodule associated to the $\big(\sur{N}_Q(T),\sur{N}_P(S)\big)$-biset $T\dom U/S$. Moreover $\Phi_{\sur{N}_Q(T)}$ is self dual, since it is the unique faithful rational irreducible representation of $\sur{N}_Q(T)$, so we can replace $\Phi_{\sur{N}_Q(T)}^*$ by $\Phi_{\sur{N}_Q(T)}$ in the right hand side of the previous isomorphism.\par
Now the biset $T\dom U/S$ splits as a disjoint union 
$$T\dom U/S=\bigsqcup_{\omega\in N_Q(T)\dom U/N_P(S)} T\dom \omega/S$$
of transitive $\big(\sur{N}_Q(T),\sur{N}_P(S)\big)$-bisets, where $N_Q(T)\dom U/N_P(S)$ is the set of $\big(N_Q(T),N_P(S)\big)$-orbits on $U$. This yields a decomposition
\begin{equation}\label{decompose}
\Q (T\dom U/S)\cong\bigoplus_{\omega\in N_Q(T)\dom U/N_P(S)}\Q \big(T\dom \omega/S\big)
\end{equation}
as $\big(\Q \sur{N}_Q(T),\Q\sur{N}_P(T)\big)$-bimodules.
\begin{mth}{Lemma}\label{libre} Let $C$ and $D$ be cyclic $p$-groups, and let $\Omega$ be a transitive $(D,C)$-biset. Then $\Phi_D\otimes_{\Q D}\Q \Omega=0$ unless $\Omega$ is left free, and $\Q\Omega\otimes_{\Q C}\Phi_C=0$ unless $\Omega$ is right free.
\end{mth}
\pf Suppose that the action of $C$ is not free. This means that $C$ is non-trivial, and that the unique subgroup $Z$ of order $p$ of $C$ acts trivially on~$\Omega$: indeed since $\Omega$ is a transitive biset, the stabilizers in $C$ of the points of~$\Omega$ are conjugate in $\mathcal{C}$, hence equal since $C$ is abelian. So these stabilizers all contain $Z$ if one of them is non trivial. Then $\Omega$ is inflated from a $(C/Z)$-set~$\sur{\Omega}$, and then $\Q\Omega\cong \Q \sur{\Omega}\otimes_{\Q(C/Z)}\Q (C/Z)$. But $\Q (C/Z)\otimes_{Q C}\Phi_C$ is the module of $Z$-coinvariants on $\Phi_C$, hence it is zero, since $\Phi_C$ is faithful. Hence $\Q\Omega\otimes_{\Q C}\Phi_C=0$ in this case. Similarly, if the action of $D$ is not free, then $\Phi_D\otimes_{\Q D}\Q\Omega=0$.\findemo
\npar\label{reduction} It follows from Lemma~\ref{libre} that to compute 
$$\Phi_{\sur{N}_Q(T)}^*\otimes_{\Q \sur{N}_Q(T)}\Q (T\dom U/S)\otimes_{\Q \sur{N}_P(S)}\Phi_{\sur{N}_P(S)}$$
using decomposition~\ref{decompose}, we can restrict to orbits $\omega=N_Q(T)uN_P(S)$, where $u\in U$, for which the $\big(\sur{N}_Q(T),\sur{N}_P(S)\big)$-biset $T\dom \omega/S$ is left and right free. The left stabilizer of the element $TuS$ of this biset is equal to
$$\{xT\in\sur{N}_G(T)\mid \exists s\in S,\;xu=us\}\mvirg$$
hence $T\dom \omega/S$ is left free if and only if
\begin{equation}\label{T puissance u}{^uS}\cap N_Q(T)\leq T\mvirg
\end{equation}
where ${^uS}=\{x\in Q\mid \exists s\in S,\;xu=us\}$ (\cite{bisetfunctors} Notation 2.3.16).\par
Similarly $T\dom \omega/S$ is right free if and only if
$$T^u\cap N_P(S)\leq S \mvirg$$
where $T^u=\{x\in P\mid\exists t\in T,\;tu=ux\}$.\par
Finally, the bimodule $L(U)_{T,S}$ is isomorphic to
$$\bigoplusc{u\in[N_Q(T))\dom U/N_P(S)]}{{^uS}\cap N_Q(T)\leq T}{T^u\cap N_P(S)\leq S}\Phi_{\sur{N}_Q(T)}\otimes_{\Q\sur{N}_Q(T)}\Q\big(T\dom N_Q(T)uN_P(S)/S\big)\otimes_{\Q \sur{N}_P(S)}\Phi_{\sur{N}_P(S)},$$
where $[N_Q(T))\dom U/N_P(S)]$ is a set of representatives of $\big(N_Q(T),N_P(S)\big)$-orbits on $U$.
\begin{mth}{Lemma} \label{cyclic} Let $p$ be an odd prime, and let $C$ and $D$ be cyclic $p$-groups. Let moreover $\Omega$ be a left and right free finite $(D,C)$-biset. Let $a\in C$, viewed as an automorphism of the $\Q C$-module $\Phi_C$. Then the image of $[\Phi_C,a]$ in $K_1(\Q D)$ by the transfer associated to the bimodule $L=\Phi_D\otimes_{\Q D} \Q\Omega$ is equal to the image of $\Ver_{\Omega}(a)\in D=D/[D,D]$ by the map $\alpha_{\Q D}$ of Remark~\ref{units}.
\end{mth}
\pf By Lemma~\ref{verlagerung} and Proposition~\ref{transfer}, we can assume that $\Omega$ is a transitive biset, of the form $(D\times C)/B$ for some subgroup $B$ of $D\times C$. Then $\Omega$ is left and right free if and only if there exists a subgroup $E$ of $C$ and an injective group homomorphism $\varphi:E\to D$ such that $B=\{\big(\varphi(e),e\big)\mid \in E\}$. There are two cases:
\begin{itemize}
\item either $E=\un$: in this case $\Omega=D\times C$, so $\Q\Omega\cong \Q D\otimes_\Q \Q C$, and $\Phi_D\otimes_{\Q D}\Q\Omega\otimes_{\Q C}\Phi_C\cong\Phi_D\otimes_\Q\Phi_C$. As a vector space over the cyclotomic field $F$ of endomorphisms of $\Phi_D$, it is isomorphic to $F\otimes_\Q\Phi_C$.  The action of $a\in C$ on this vector space is given by the matrix of $a$ acting on $\Phi_C$.\par
Suppose that $a$ is a generator of $C$, of order $p^n$. Then this action is the action by multiplication of a primitive $p^n$-th root of unity $\zeta$ on the field $\Q(\zeta)$. As an element of $K_1(F)$, it is equal to the determinant of the matrix representing this multiplication, i.e. to the norm $N_{\Q(\zeta)/\Q}(\zeta)$, which is equal to 1, as the $p^n$-th cyclotomic polynomial has even degree $p^{n-1}(p-1)$ and value 1 at 0. It follows that $[\Phi_C,a]$ is mapped to the identity element of $K_1(\Q D)$ in this case. Since this holds for a generator $a$ of $C$, the same is true for any element $a$ of $C$.\par
In this case also, a set of representatives of $[D\dom \Omega]$ is the set ${1}\times C$, which is invariant by right multiplication by $C$. It follows that the elements $d_{a,x}\in D$ defined for $a\in C$ and $x\in[D\dom\Omega]$ by $xa=d_{a,x}x'$, for $x'\in [D\dom\Omega]$, are all equal to 1. So the transfer $\Ver_\Omega$ is also the trivial homomorphism in this case.
\item or $E\neq\un$: let $Z$ denote the unique subgroup of order $p$ of $C$. Tensoring over $\Q C$ the exact sequence of $(\Q C,\Q C)$-bimodules
$$0\to \Phi_C\to\Q C\to \Q (C/Z)\to 0$$
with $\Phi_D\otimes_{\Q D}\Q \Omega$ gives the exact sequence
$$0\to \Phi_D\otimes_{\Q D}\Q \Omega\otimes_{\Q C}\Phi_C\to \Phi_D\otimes_{\Q D}\Q \Omega\to \Phi_D\otimes_{\Q D}\Q (\Omega/Z)\to 0\mpoint$$
But $\Omega/Z$ is not free as a left $D$-set, since the unique subgroup $\varphi(Z)$ of order $p$ of $D$ stabilizes $BZ\in \Omega/Z$, as for $z\in Z$ and $e\in E$
$$\varphi(z)\big(\varphi(e),e\big)=\big(\varphi(ze),e\big)=\big(\varphi(ze),ze\big)z\mpoint$$
By Lemma~\ref{libre}, it follows that $\Phi_D\otimes_{\Q D}\Q (\Omega/Z)=0$, hence
$$\Phi_D\otimes_{\Q D}\Q \Omega\otimes_{\Q C}\Phi_C\cong \Phi_D\otimes_{\Q D}\Q \Omega\mpoint$$
As a vector space over the cyclotomic field $F$ of endomorphisms of $\Phi_D$, this is isomorphic to $F\otimes_\Q\Q[D\dom\Omega]$. The action of $a\in C$ on this vector space is given for $x\in [D\dom\Omega]$ and $\lambda\in F$ by
$$(\lambda\otimes x)a=\lambda\otimes xa=\lambda\otimes d_{a,x}\sigma_a(x)=\lambda d_{a,x}\otimes \sigma_a(x)\mvirg$$
where $d_{a,x}\in D$ and $\sigma_a(x)\in [D\dom\Omega]$ are defined by $xa=d_{a,x}\sigma_a(x)$. In other words, the matrix of the action of $a$ is the product of the permutation matrix of $\sigma_a$ with a diagonal matrix of coefficients $d_{a,x}$, for $x\in [D\dom\Omega]$. In $K_1(F)$, this matrix is equal to its determinant, that is the signature of $\sigma_a$, which is equal to 1 as $\sigma_a$ is a product of cycles of odd length (equal to some power of $p$), multiplied by the product of the elements $d_{a,x}$, that is the image in $K_1(\Q D)$ of $\Ver_\Omega(a)$, as was to be shown.\findemo 
\end{itemize}
\begin{rem}{Remark} Recall that if $Q$ is a central subgroup of finite index $n$ in a group $G$, then the transfer $G/[G,G]\to Q$ is induced by the map $g\mapsto g^n$ from $G$ to $Q$ (see \cite{rotman} Theorem 7.47). It follows easily that in the situation of Lemma~\ref{cyclic}, if $\Omega=(D\times C)/B$, where $B=\{\big(\varphi(e),e\big)\mid e\in E\}$ for a subgroup $E$ of $C$ and an injective homomorphism $\varphi:E\to D$, the transfer $\Ver_\Omega:C\to D$ is given by $a\mapsto \varphi(a^{|C:E|})$. Moreover $|C:E|=|D\dom\Omega|$.
\end{rem}
\begin{mth}{Theorem} \label{formula} Let $p$ be an odd prime. Let $P$ and $Q$ be finite $p$-groups, and let $U$ be a finite $(Q,P)$-biset. 
\begin{enumerate}
\item Let $S$ be a genetic subgroup of $P$ and $T$ be a genetic subgroup of $Q$. Let $\mathcal{D}(U)_{T,S}$ be the set of orbits $N_Q(T)uN_P(S)$ of those $u\in U$  for which {$T^u\cap N_P(S)\leq S$} and $^uS\cap N_Q(T)\leq T$ (see~\ref{T puissance u} for notation). Then for $\omega\in\mathcal{D}(U)_{T,S}$, the set $T\dom\omega/S$ is a left and right free $\big(N_Q(T)/T,N_P(S)/S\big)$-biset, and the map
$$\Gamma(U)_{T,S}:N_P(S)/S\to N_Q(T)/T$$ 
sending $a\in N_P(S)/S$ to
$$\prod_{\omega\in\mathcal{D}(U)_{T,S}}\Ver_{T\dom\omega/S}(a)$$
is a well defined group homomorphism.
\item Let $\mathcal{B}_P$ and $\mathcal{B}_Q$ be genetic bases of $P$ and $Q$, respectively. Then the map $\Gamma(U):\Gamma(P)\to\Gamma(Q)$ giving the biset functor structure of $\Gamma$ is the map
$$\Gamma(P)=\prod_{S\in \mathcal{B}_P}\big(N_P(S)/S\big)\to\prod_{T\in \mathcal{B}_Q}\big(N_Q(T)/T\big)=\Gamma(Q)$$
with component $(T,S)$ equal to $\Gamma(U)_{T,S}$.
\end{enumerate}
\end{mth}
\pf This results from Paragraph~\ref{functor}, Lemma~\ref{libre}, Paragraph~\ref{reduction}, and Lemma~\ref{cyclic}.\findemo
\section{Examples}
\begin{mth}{Proposition} Let $P$ be a finite $p$-group, for $p$ odd, and let $\mathcal{B}$ be a genetic basis of~$P$. Let $N\normal P$, and $\sur{P}=P/N$.
Let $\mathcal{B}_N$ be the subset of $\mathcal{B}$ defined by
$$\mathcal{B}_N=\{S\in\CB\mid S\geq N\}\mpoint$$
Then:
\begin{enumerate}
\item The set $\sur{\CB}=\{\sur{S}=S/N\mid S\in \CB_N\}$ is a genetic basis of $\sur{P}$.
\item Up to the identification of $N_{\sur{P}}(\sur{S})$ with $N_P(S)/S$, for $S\in\CB_N$, the inflation morphism
$$\Inf_{P/N}^P:\Gamma(P/N)=\prod_{\sur{S}\in\sur{\CB}}\limits\big(N_{\sur{P}}(\sur{S}\big)/\sur{S})\to \Gamma(P)=\prod_{S\in\CB}\limits\big(N_P(S)/S\big)$$
is the embedding in the product of the factors of $\Gamma(P)$ corresponding to genetic subgroups $S$ containing $N$.
\item Similarly, the deflation morphism
$$\Def_{P/N}^P:\Gamma(P)=\prod_{S\in\CB}\limits\big(N_P(S)/S\big)\to \Gamma(P/N)=\prod_{\sur{S}\in\sur{\CB}}\limits\big(N_{\sur{P}}(\sur{S}\big)/\sur{S})$$
is the projection onto the product of the factors of $\Gamma(P)$ corresponding to genetic subgroups $S$ containing $N$.
\end{enumerate}
\end{mth}
\pf Assertion 1 is clear from the definitions: if $N\leq S\leq P$, then $S$ is genetic in $P$ if and only if $S/N$ is genetic in $P/N$. Moreover the relation $\bizlie{P/N}$ gives the relation $\bizlie{P}$ by inflation.\par
Now the inflation morphism $\Inf_{P/N}^P$ is defined by the $(P,\sur{P})$-biset $U=P/N$, with natural actions of $P$ and $\sur{P}$. Let $u=xN$ be an element of $U$, for some $x\in P$. Let $T\in\CB$ and $\sur{S}\in\sur{\CB}$.
\begin{eqnarray*}
T^u\cap N_{\sur{P}}(\sur{S})&=&\{\sur{y}=yN\in N_{\sur{P}}(\sur{S})\mid\exists t\in T,\;tu=u\sur{y}\}\\
&=&\{\sur{y}=yN\in N_{\sur{P}}(\sur{S})\mid\exists t\in T,\;txN=xyN\}\\
&=&\{\sur{y}=yN\in N_{\sur{P}}(\sur{S})\mid yN\in {T^x}\,N\}\\
&=&\big({T^x}\,N\cap N_P(S)\big)/N=\big(T^x\cap N_P(S)\big)N/N\mvirg
\end{eqnarray*}
where the last two equalities hold because $S\geq N$. Hence $T^u\cap N_{\sur{P}}(\sur{S})\leq\sur{S}$ if and only if $T^x\cap N_P(S)\leq S$.\par
On the other hand
\begin{eqnarray*}
{^u\sur{S}}\cap N_P(T)&=&\{y\in N_P(T)\mid \exists \sur{s}=sN\in\sur{S},\;yu=u\sur{s}\}\\
&=&\{y\in N_P(T)\mid \exists s\in S,\;yxN=xNs\}\\
&=&\{y\in N_P(T)\mid y\in {^xS}\}\\
&=&{^xS}\cap N_P(T)\mpoint
\end{eqnarray*}
Hence ${^u\sur{S}}\cap N_P(T)\subseteq T$ if and only if ${^xS}\cap N_P(T)\leq T$. \par
If moreover $T^u\cap N_{\sur{P}}(\sur{S})\leq\sur{S}$, i.e. $T^x\cap N_P(S)\leq S$, it follows that $T\bizlie{P}S$, hence $T=S$ since $T$ and $S$ belong to the same genetic basis $\CB$. Moreover $x\in N_P(S)$, and the induced group homomorphism $N_{\sur{P}}(\sur{S})/\sur{S}\to N_P(T)/T$ is the canonical isomorphism $N_{\sur{P}}(\sur{S})/\sur{S}\to N_P(S)/S$. This completes the proof of Assertion~2.\par
For Assertion 3, we consider the deflation map $\Def_{P/N}^P:\Gamma(P)\to\Gamma(P/N)$. It corresponds to the biset $V=P/N$, with left action of $\sur{P}$ and right action of $P$. For $v=yN\in V$, for $T\in\CB$ and $\sur{S}\in\sur{\CB}$, and with the same notation as above, the computation is similar: we have $T^v\cap N_{\sur{P}}(\sur{S})\leq\sur{S}$ if and only if $T^y\cap N_P(S)\leq S$, and ${^v\sur{S}}\cap N_P(T)\leq T$ if and only if ${^yS}\cap N_P(T)\leq T$. These two conditions are fulfilled if and only if $S=T$ and $y\in N_P(S)$. This completes the proof.\findemo\vspace{2ex}
Recall that the faithful part $\partial F(P)$ of the evaluation of a biset functor~$F$ at a group $P$ is the set of faithful elements of $F(P)$, introduced in~\cite{bisetfunctors} Definition~6.3.1: it is the set of elements $u\in F(P)$ such that $\Def_{P/N}^Pu=0$ for any non trivial normal subgroup $N$ of $P$. Equivalently $\Def_{P/Z}^Pu=0$ for any non trivial central subgroup $Z$ of $P$. The following is now clear:\pagebreak[3]
\begin{mth}{Corollary}\label{Gamma faithful} Let $p$ be an odd prime and $P$ be a finite $p$-group. Let $\CB$ be a genetic basis of $P$. Then the faithful part $\partial\Gamma(P)$ of $\Gamma(P)$ is equal to
$$\partial \Gamma(P)=\prodb{S\in \CB}{S\cap Z(P)=\un}\big(N_P(S)/S\big)\mpoint$$
\end{mth}
\section{Genome and Roquette category}\label{roquette}
\npar Let $F$ be a $p$-biset functor. It is shown in \cite{bisetfunctors} Theorem 10.1.1 that if $P$ is a finite $p$-group and $\CB$ is a genetic basis of $P$, then the map
$$\mathcal{I}_\CB=\bigoplus_{S\in\CB}\Indinf_{N_P(S)/S}^P:\bigoplus_{S\in \CB}\partial F\big(N_P(S)/S\big)\to F(P)$$
is always split injective. When $\mathcal{I}_\CB$ is an isomorphism for one particular genetic basis $B$ of $P$, then $\mathcal{I}_{\CB'}$ is an isomorphism for any other genetic basis $\CB'$ of $P$.\par
The functors for which $\mathcal{I}_B$ is an isomorphism for any finite $p$-group $P$ and any genetic basis $\CB$ of $P$ are called {\em rational} $p$-biset functors. It has been shown further (\cite{roquette-category}) that these rational $p$-biset functors are exactly those $p$-biset functors which factorize through {\em the Roquette category} $\mathcal{R}_p$ of $p$-groups: more precisely (\cite{roquette-category}, Definition 3.3), the category $\mathcal{R}_p$ is defined as the idempotent additive completion of a specific quotient $\mathcal{R}_p^\sharp$ of the category $\mathcal{C}_p$, so there is a canonical additive functor $\pi_p:\mathcal{C}_p\to\mathcal{R}_p$, equal to the composition of the projection functor $\mathcal{C}_p\to\mathcal{R}_p^\sharp$ and the inclusion functor $\mathcal{R}_p^\sharp\to \mathcal{R}_p$. The rational $p$-biset functors are the additive functors $F:\mathcal{C}_p\to\Ab$ for which there exists an additive functor $\sur{F}:\mathcal{R}_p\to\Ab$ such that $F=\sur{F}\circ\pi_p$. In this case, the functor~$\sur{F}$ is unique.
\begin{mth}{Proposition} Let $p$ be an odd prime. Then the genome $p$-biset functor $\Gamma$ is rational. 
\end{mth}
\pf Let $P$ be a $p$-group, and $\CB$ be a genetic basis of $P$. If $S\in \CB$, then $Q=N_P(S)/S$ is cyclic, so the trivial subgroup $S/S$ of $Q$ is the only one intersecting trivially the center of $Q$, and it is a genetic subgroup of $Q$. By Corollary~\ref{Gamma faithful}, we have that 
$$\partial\Gamma\big(N_P(S)/S\big)=N_P(S)/S\mpoint$$
Now the induction-inflation map $\Indinf_{N_P(S)/S}^P$ is given by the $\big(P,N_P(S)/S\big)$-biset $U=P/S$. Let $T\in\CB$, and let $u=xS\in U$ such that
$$T^u\cap N_Q(S/S)\leq S/S\;\;\; \hbox{and}\;\;\;{^u(S/S)}\cap N_P(T)\leq T\mpoint$$
The first inclusion means that
$$\{yS\in N_P(S)/S\mid\exists t\in T,\;txS=xSy\}=S/S\mpoint$$
In other words $T^x\cap N_P(S)\leq S$. The second inclusion means similarly that
$$\{y\in N_P(T)\mid \exists s\in S/S,\;txS=xSs\}\leq T\mvirg$$
that is $N_P(T)\cap {^xS}\leq T$. Hence $T\bizlie{P}S$, thus $T=S$ since $T$ and $S$ both belong to a genetic basis of $P$. Moreover $x\in N_P(S)$, and the morphism we get from $N_P(S)/S$ to $N_P(T)/T$ is the identity map. \par
In other words, the map 
$$\Indinf_{N_P(S)/S}^P:\partial\Gamma\big(N_P(S)/S\big)=N_P(S)/S\to \Gamma(P)$$
is the canonical embedding of $N_P(S)/S$ in $\Gamma(P)$. It clearly follows that the map $\mathcal{I}_B$ is an isomorphism, hence $\Gamma$ is rational.\findemo
\begin{mth}{Corollary} Let $p$ be an odd prime. Then there exists a unique additive functor $\sur{\Gamma}$ from the Roquette category $\mathcal{R}_p$ to $\Ab$ such that $\Gamma=\sur{\Gamma}\circ \pi_p$. Moreover $\sur{\Gamma}(\partial P)=\partial\Gamma(P)$ for any finite $p$-group $P$, where $\partial P$ is the edge of~$P$ in $\mathcal{R}_p$. In particular $\sur{\Gamma}(\partial C)=C$ for any cyclic $p$-group $C$.
\end{mth}
\pf This follows from the definition and properties of the category $\mathcal{R}_p$ (for the definition of the {\em edge} $ \partial P$ of a $p$-group $P$ in the Roquette category, see \cite{roquette-category} Definition~3.7).\findemo

\begin{thebibliography}{10}

\bibitem{barker-genotype}
L.~Barker.
\newblock Genotypes of irreducible representations of finite $p$-groups.
\newblock {\em Journal of Algebra}, 306:655--681, 2007.

\bibitem{bisetfunctors}
S.~Bouc.
\newblock {\em Biset functors for finite groups}, volume 1990 of {\em Lecture
  Notes in Mathematics}.
\newblock Springer, 2010.

\bibitem{roquette-category}
S.~Bouc.
\newblock The {R}oquette category of finite $p$-groups.
\newblock {\em Journal of the European Mathematical Society}, 17:2843--2886,
  2015.

\bibitem{bouc-romero}
S.~Bouc and N.~Romero.
\newblock The {W}hitehead group of (almost) extra-special $p$-groups with
  $p$-odd.
\newblock Preprint, arXiv:1604.06306, 2016.

\bibitem{curtisreinerII}
C.~Curtis and I.~Reiner.
\newblock {\em Methods of representation theory with applications to finite
  groups and orders}, volume \MakeUppercase{\romannumeral2} of {\em Wiley
  classics library}.
\newblock Wiley, 1990.

\bibitem{curtisreiner}
C.~Curtis and I.~Reiner.
\newblock {\em Methods of representation theory with applications to finite
  groups and orders}, volume \MakeUppercase{\romannumeral1} of {\em Wiley
  classics library}.
\newblock Wiley, 1990.

\bibitem{bob}
R.~Oliver.
\newblock {\em Whitehead groups of finite groups}, volume 132 of {\em London
  Mathematical Society Lecture Note Series}.
\newblock Cambridge University Press, Cambridge, 1988.

\bibitem{romero-whitehead}
N.~Romero.
\newblock Computing {W}hitehead groups using genetic bases.
\newblock {\em Journal of Algebra}, 450:646--666, 2016.

\bibitem{roquette}
P.~Roquette.
\newblock Realisierung von {D}arstellungen endlicher nilpotenter {G}ruppen.
\newblock {\em Arch. Math.}, 9:224--250, 1958.

\bibitem{rotman}
J.~J. Rotman.
\newblock {\em An introduction to the theory of groups}, volume 148 of {\em
  Graduate Texts in Mathematics}.
\newblock Springer-Verlag, New York, fourth edition, 1995.

\end{thebibliography}

Serge Bouc\\
LAMFA-CNRS UMR7352\\
33 rue St Leu, Amiens Cedex 01\\
France\\
{\tt email: serge.bouc@u-picardie.fr}
\end{document}